\documentclass[11pt]{article}

\input{amssym.def}
\input{amssym}
\usepackage{eufrak}

\newtheorem{theorem}{Theorem}
\newtheorem{corollary}[theorem]{Corollary}
\newtheorem{definition}[theorem]{Definition}

\newtheorem{lemma}[theorem]{Lemma}
\newtheorem{proposition}[theorem]{Proposition}
\newtheorem{remark}[theorem]{Remark}

\def\hl{\hat{l}}
\def\hL{\hat{L}}

\def\un{\underline}
\def\qq{q^{-1}}

\def\Tr{\mathrm{Tr}}

\def\hLL{\hat{{\cal L}}}
\def\LL{{\cal{L}}}

\def\x{\mathsf{x}}
\def\y{\mathsf{y}}

\def\de{\delta}

\def\ot{\otimes}
\def\C{{\Bbb C}}

\def\vv{V^{\otimes 2}}

\def\om{\omega}

\def\ov{\overline}
\def\un{\underline}
\def\hmu{\hat{\mu}}
\def\la{{\lambda}}
\def\H{{\cal{H}}}

\def\be{\begin{equation}}
\def\ee{\end{equation}}
\def\hf{\hfill \rule{6.5pt}{6.5pt}}

\begin{document}

\makeatletter
\renewcommand{\theequation}{{\thesection}.{\arabic{equation}}}
\@addtoreset{equation}{section} \makeatother

\title{Representations of Spectrum of $GL(m)$ type Quantum Matrices}
\author{\rule{0pt}{7mm} Dimitry
Gurevich\thanks{gurevich@ihes.fr}\\
{\small\it Higher School of Modern Mathematics,}\\
{\small\it Moscow Institute of Physics and Technology}\\
\rule{0pt}{7mm} Pavel Saponov\thanks{Pavel.Saponov@ihep.ru}\\
{\small\it
National Research University Higher School of Economics,}\\
{\small\it 20 Myasnitskaya Ulitsa, Moscow 101000, Russian Federation}\\
{\small \it and}\\
{\small \it
Institute for High Energy Physics, NRC "Kurchatov Institute"}\\
{\small \it Protvino 142281, Russian Federation}\\
\rule{0pt}{7mm} Mikhail Zaitsev\thanks{mzaicev02@gmail.com}\\
{\small\it National Research University Higher School of Economics,}\\
{\small\it 20 Myasnitskaya Ulitsa, Moscow 101000, Russian Federation}}

\maketitle

\begin{abstract}
In the present paper we are dealing with reflection equation algebras $\LL(R)$ corresponding to even skew-invertible Hecke 
symmetries. Our main result consists in computing the characters of the spectral values of the generating matrix $L$ of $\LL(R)$
in finite-dimensional representations labeled by partitions of integers. As is known, the spectral values belong to an algebraic 
extension of the center of the reflection equation algebra and elements of the center can be presented as 
symmetric functions in spectral values. As an application of our approach, we calculate the characters of the
power sums $\mathrm{Tr}_R(L^n)$ in the mentioned finite dimensional representations.

In a particular case of the Drinfeld-Jimbo $R$-matrix the enveloping algebra $U(gl(N))$ can be obtained as a specific limit
of the reflection equation algebra. In this limit our results for power sums coincide with those obtained in \cite{PP}.
\end{abstract}

{\bf AMS Mathematics Subject Classification, 2020:} 81R60


\section{Introduction}
By a quantum matrix we mean the matrix $L$ composed of generators of an algebra, associated with two compatible braidings 
$R$ and $F$ (see \cite{IOP} for definitions and details). In the current paper we are dealing with a particular subclass of such algebras,
corresponding to the case $F=R$. Such algebras are called Reflection Equation (RE)
algebras and denoted $\LL(R)$.  A common feature of each such an algebra is that its generating matrix $L$ meets a polynomial relation
with coefficients which are central in $\LL(R)$. This relation is a noncommutative analog of the classical Cayley-Hamilton
identity. This property enables us to introduce eigenvalues of the generating matrix $L$ as elements of the central extension of the algebra
$\LL(R)$. The aim of the current paper is to compute the characters of the eigenvalues of the matrix $L$ in a family of (conjecturally
irreducible) representations of the algebra  $\LL(R)$ for the case when $R$ is an even Hecke symmetry of the bi-rank $(m|0)$ 
(see Section 2). The algebra $\LL(R)$ corresponding to such a Hecke symmetriy is called the RE algebra of $GL(m)$ type.

We illustrate the above constructions on the example of the algebra $U(gl(N))$, which can be obtained as a specific limit
of the RE algebra, associated with $R$ coming from the $U_q(sl(N))$ --- the so-called Drinfeld-Jimbo $R$-matrix 
(see Section 3).

Let $\hl_i^j$, $1\leq i, j \leq N$ be the standard generators of $U(gl(N))$ satisfying the system of commutation relations:
$$
\hl_i^j\, \hl_k^m-l_k^m\, \hl_i^j=\hl_i^m\, \de_k^j- \hl_k^j\, \de_i^m.
$$
This system can be written in the following matrix form:
\be
P\, (\hL\ot I)\, P\, (\hL\ot I)-(\hL\ot I)\, P\, (\hL\ot I)\, P=P\, (\hL\ot I)- (\hL\ot I)\, P,
\label{odin}
\ee
where $\hL=\|\hl_i^j\|_{1\leq i,j \leq N}$, $I$ is the $N\times N$ identity matrix, and $P$ is the $N^2\times N^2$ matrix of the usual flip.
The matrix $\hL$ will be called {\it the generating matrix} of the algebra $U(gl(N))$.

The generating matrix $\hL$  is subject to a polynomial relation which is an analog of the usual Cayley-Hamilton identity:
\begin{equation}
\hL^N +\sum_{k=1}^N(-1)^ka_k(\hL)\,\hL^{N-k} = 0,
\label{CH}
\end{equation}
where the coefficients $a_k(\hL)$ belong to the center $Z(U(gl(N))$ of the algebra $U(gl(N))$.

We introduce the {\it spectral values} (or eigenvalues) $\hmu_i$, $1\leq i\leq N$, of the generating matrix $\hL$ in the natural way:
\be
a_1(\hL) = \sum_{i=1}^N\, \hmu_i,
\qquad a_2(\hL) = \sum_{i<j}^N\, \hmu_i\, \hmu_j, \qquad \dots\qquad  a_N(\hL) = \prod_{i=1}^N\hmu_i.
\label{sta}
\ee
The eigenvalues $\hat \mu_i$ are assumed to be central in the extended algebra $U(gl(N))[\hmu_1,\dots,\hmu_N]$.
Due to definition (\ref{sta}) the Cayley-Hamilton identity (\ref{CH}) can be written in a factorized form in the extended algebra:
$$
\prod_{k=1}^N(\hL - \hmu_k\,I) = 0.
$$
Moreover, any central element of the algebra $U(gl(N))[\hmu_1,\dots,\hmu_N]$ can be expressed as a symmetric function in the
eigenvalues $\hmu_k$.

It is known that any finite dimensional irreducible representation of the algebra $U(gl(N))$ (considered up to isomorphisms) can be
labeled by a partition 
$$
\la=(\la_1\geq \la_2\geq \dots\geq \la_N\ge 0).
$$
Since the algebra $gl(N)$ is reductive,
the image $\pi_\la(z)$ of any central element $z\in  Z(U(gl(N)))$ in such a representation $\pi_\la$ is a scalar operator. The values
of the images $\pi_\la(p_k)$ for the central elements $p_k=\mathrm{Tr}\, \hL^k$, $k\ge 0$, which are called {\em power sums}, were
computed in \cite{PP}.

It is not difficult to  assign a value $\pi_\la(\hmu_i)$  to each eigenvalue  $\hmu_i$, which would be compatible  with the values
$\pi_\la(p_k )$ or, what is equivalent, compatible with the values of $\pi_\la(a_k)$. In this sense we speak about representations
of the spectrum of the matrix $\hL$.

Let us pass now to a general  case. By a modified RE algebra we mean an algebra $\hLL(R)$, such that its generating matrix $\hL$ 
meets the relation  analogous to (\ref{odin}), where the matrix $P$ is replaced by a Hecke symmetry $R$. If the matrix $R$ is a 
deformation of the usual flip (for instance $R$ comes from the Quantum Group $U_q(sl(N))$ or it is a Crammer-Gervais symmetry),
the algebra $\hLL(R)$ is a deformation of the algebra $U(gl(N))$.  

Note that for a Hecke symmetry $R$ the modified RE algebra $\hLL(R)$ is isomorphic to the RE algebra $\LL(R)$ defined by purely
quadratic relations on generators (see Section 3 for detail).  

Another important property of the $GL(m)$ type RE algebras $\hLL(R)$ is that their representation theory is similar to that of 
$U(gl(m))$. The corresponding finite dimensional representations $\pi_\lambda$ are  labeled by partitions of integers with height
(the number of nonzero components) no greater than $m$. 
Note, that the number $m$ is determined by the bi-rank of $R$ (see Section 2). 

The paper is organized as follows. In Section 2 we recall some definitions and properties of the Hecke algebras and their $R$-matrix
representations. Section 3 is devoted to a review of certain properties of the centers of the RE algebras and of
some their finite dimensional representations. The last section  contains the main results of the paper: we compute  the characters
of the spectral elements of the quantum matrices $L$ and $\hL$ in the mentioned representations  of the RE algebras and show
the limit transition to the classical results of \cite{PP}. Moreover, for the matrix $R$ coming from the quantum group $U_q(sl(N))$
(in this case $m=N$) our result for the characters of power sums $\mathrm{Tr}_R(L^k)$ coincides with that obtained in \cite{JLM}
by different methods.

\section{Hecke algebras and Hecke symmetries}

In this section we recall some facts on Hecke algebras and Hecke type braidings, called Hecke symmetries, which are essential for
studying the center of RE algebra $\LL(R)$ and its finite dimensional representations.

The Hecke algebra ${\cal H}_n(q)$ is an associative unital algebra over the complex field ${\Bbb C}$ generated by the finite set of
Artin's generators $\tau_i$, $1\le i\le n-1$, subject to the relations
$$
\tau_i\, \tau_{i+1}\, \tau_i=\tau_{i+1}\, \tau_{i}\, \tau_{i+1},\qquad \tau_i\, \tau_j=\tau_j\, \tau_i\quad |i-j|\geq 2,
$$
and
$$
(q\, e-\tau_i)(\qq\, e+\tau_i)=0,\quad 1\leq i \leq n-1,\qquad q\in {\Bbb C}\setminus 0.
$$
Here $e$ is the unit element of the algebra.
Note that for $q=\pm1$ we get the defining relations of the symmetric group algebra ${\Bbb C}[S_n]$.

In what follows we assume the parameter $q$ to be generic which means that $q^{k}\not = 1$, $\forall  k\ge 1$. For a generic $q$
the Hecke algebra ${\cal H}_n(q)$ is isomorphic to ${\Bbb C}[S_n]$ for all $n\ge 1$. Therefore, the Hecke algebra ${\cal H}_n(q)$
is finite dimensional, semisimple and isomorphic to a direct product of finite dimensional matrix algebras. As a consequence, in the
algebra ${\cal H}_n(q)$ there exists a linear basis of  ``matrix units'' $\{e^\lambda_{ij}\}_{1\le i,j\le d_\lambda}$ parameterized by
all partitions $\lambda$ of the integer $n$. Here $d_\lambda$ is the number of standard Young tables corresponding to the
partition $\lambda$. Note that the multiplication table for these matrix units  reads:
\be
e^\la_{ki}\,e^\mu_{r p}=\delta^{\la\mu}\delta_{ir}\,e^\la_{kp}.
\label{prod-mun}
\ee
The elements $e_i^\la:=e^\la_{ii}$, $\lambda\vdash n$, $1\le i \le d_\la$, are {\em the primitive orthogonal idempotents} defining
a resolution of the unit $e$:
\be
e=\sum_{\lambda\,\vdash n} \sum_{i=1}^{d_\lambda}e^\la_i.
\label{reso}
\ee
An explicit construction of all these  idempotents is exhibited in \cite{OP} (see also the references therein).

A maximal commutative subalgebra of ${\cal H}_n(q)$ is generated by the so-called {\it Jucys-Murphy elements} $\{j_k\}_{1\le k\le n}$,
which are defined by the relations:
\be
j_1 = e,\qquad j_k = \tau_{k-1}\,j_{k-1}\tau_{k-1},\quad 2\le k\le n.
\label{JM}
\ee

Let $V$ be a finite dimensional complex vector space
$\dim_{\,\Bbb C} V=N$. A {\it braiding} is an operator $R:\vv\to \vv$ meeting the {\it braid relation}:
\be
R_{12}R_{23}R_{12}= R_{23}R_{12}R_{23}.
\label{braid}
\ee
Here $R_{12} = R\otimes \mathrm{Id}_V$, $R_{23} = \mathrm{Id}_V\otimes R$.

Hereafter, we fix a basis $\{x_i\}_{1\le i\le N} $ in the space $V$ and we identify operators in spaces
$V^{\otimes p}$, $p\ge 1$, with their matrices in the corresponding tensor basis $\{x_{i_1}\ot\dots \otimes x_{i_p}\}$.

Below we use the following embeddings of the $N^2\times N^2$ matrix $R$ into the space of $N^p\times N^p$ matrices:
\be
R_{i\,i+1} =  I^{\otimes (i-1)}\otimes R\otimes I^{\otimes (p-i-1)},\qquad 1\le i\le p-1,
\label{R-emb}
\ee
$I$ being the $N\times N$ identity matrix. In what follows we do not fix the integer $p$ just assuming it to be sufficiently large in order
all the matrix formulas make sense. Besides, we shall use a shorthand  notation $R_i$ for $R_{ii+1}$.
As can be easily seen the braid relation (\ref{braid}) takes place for all copies $R_i$:
\be
R_iR_{i+1}R_i = R_{i+1}R_iR_{i+1},\quad \forall\, i\ge 1.
\label{braid-i}
\ee

A braiding $R$ meeting the Hecke condition
\be
(q\, I^{\ot 2}-R)(\qq I^{\ot 2}+R)=0,\qquad q\in \C\setminus 0,\quad q\not=\pm 1
\label{above}
\ee
is called {\it a Hecke symmetry}. If $q=\pm 1$, the corresponding braiding $R$ is called {\it an involutive symmetry}.

Any Hecke symmetry $R$ defines an $R$-matrix representation $\rho_R: {\cal H}_n(q)\rightarrow \mathrm{End}(V^{\ot n})$
of the Hecke algebras. On Artin's generators this representation is given by the map:
\be
\rho_R(\tau_k)\mapsto R_k: = I^{\otimes (k-1)}\otimes  R\otimes I^{\otimes (n -k-1)},\quad 1\le k\le n-1.
\label{R-rep}
\ee
The images of the Jucys-Murphy elements in the $R$-matrix representation are given by the following matrices:
\be
J_k = \rho_R(j_k) = R_{k-1}\dots R_2R_1^2R_2\dots R_{k-1}.
\label{JM-rep}
\ee
Below the matrices $J_k$ will be also referred to as Jucys-Murphy elements.

With any symmetry $R$ one can naturally associate the $R$-symmetric and  $R$-skew-symmetric algebras of the space $V$
(see \cite{G,GPS1} for definitions). A symmetry $R$ is called {\em even} if the $R$-skew-symmetric algebra is finite-dimensional.
Then the Hilbert-Poincare series of the $R$-skew-symmetric algebra is a polynomial of a degree $m$. In this case 
the symmetry $R$ is called to be of the bi-rank $(m|0)$ (see \cite{GPS1} for detail).

Moreover, we  assume any symmetry $R$ to be skew-invertible. This means that there exists an operator
$\Psi:V^{\otimes 2}\to V^{\otimes 2}$
such that
$$
\mathrm{Tr}_{(2)}( R_{12} \Psi_{23})=P_{13} = \mathrm{Tr}_{(2)}(\Psi_{12}R_{23}).
$$

Below we also use the $N\times N$ matrices $B=\|B_i^j\|$ and $C = \|C_i^j\|$, defined by
\be
B_i^j = \sum_{a=1}^N \Psi_{ai}^{aj},\qquad C_i^j = \sum_{a=1}^N \Psi_{ia}^{ja},
\label{C}
\ee
where $\|\Psi_{ij}^{ks}\|$ is a matrix of the operator $\Psi$ in the tensor basis of the space $V^{\otimes 2}$.

\begin{definition}\rm
Let $M$ be an $N\times N$ matrix with entries form any complex vector space $U$. Its $R$-trace is defined as
the map $\mathrm{Mat}_N(U)\rightarrow U$ given by the rule:
\be
{\Tr_R}\, M:=\Tr\,( C\cdot M).
\label{tra}
\ee
\end{definition}

Below we use the convenient shorthand notation for the multiple trace:
$$
\langle M_{12\dots k}\rangle_{1\dots k}:=
\Tr_{R(12\dots k)}M_{12\dots k} = \Tr_{(12\dots k)}(C_1C_2\dots C_kM_{12\dots k})
$$
where  $M_{12\dots k}\in \mathrm{Mat}_N(U)^{\otimes k}$ is an arbitrary $N^k\times N^k$ matrix and
$C_i = I^{\otimes(i-1)}\otimes C\otimes I^{\otimes(k-i)}$.

The $R$-trace possesses the useful {\it cyclic property}. Namely, for any polynomial $Q$ in matrices $R_1,\dots, R_{k-1}$ and
for an arbitrary matrix $M_{12\dots k}$
the following identity takes place:
$$
\langle M_{12\dots k} \,Q(R_1,\dots, R_{k-1})\rangle_{1\dots k} = \langle Q(R_1,\dots, R_{k-1})\,M_{12\dots k}\rangle_{1\dots k}.
$$
This cyclic property is an immediate consequence of the following relation:
$$
R_i C_iC_{i+1} = C_iC_{i+1} R_i, \quad \forall\, i\ge 1.
$$

\section{Reflection equation algebra and its basic properties}
\label{sec:3}

Now, we define a Reflection Equation (RE) algebra associated with a skew-invertible Hecke symmetry $R$. By definition, an RE
algebra $\LL(R)$, corresponding to a given symmetry $R$, is a unital associative algebra finitely generated by the elements
$l_i^j$, $1\le i,j\le N$, subject to the quadratic permutation relations:
\be
R\, (L\ot I)R\, (L\ot I)-(L\ot I) R\, (L\ot I) R=0,\qquad L = \|l_i^j\|.
\label{odinn}
\ee
Note that the generating matrix $L= \sum_{i,j = 1}^N E_i^j\otimes l_i^j$ is actually an element of the algebra $\mathrm{End}({\Bbb
C}^N)\otimes \LL(R)$, where $E_i^j$ are the matrix units -- elements of a linear basis of $\mathrm{End}({\Bbb C}^N)$.

Let us introduce a new generating matrix $\hL$:
\be
L=Ie_{\cal L}- (q-\qq) \hL,
\label{isom}
\ee
where $e_{\cal L}$ is the unit element of the algebra $\LL(R)$.
As a consequence of (\ref{odinn}) and of the Hecke condition (\ref{above}) we get the following quadratic-linear relation for the
matrix $\hL$:
\be
R\, (\hL\ot I) R\, (\hL\ot I)-(\hL\ot I) R\, (\hL\ot I) R = R\, (\hL\ot I)- (\hL\ot I)  R,
\label{mREA}
\ee
which is similar to (\ref{odin}), but with $P$ replaced by $R$. The corresponding algebra is called a {\em  modified} RE algebra
and is denoted $\hLL(R)$.
Evidently, the algebras $\LL(R)$ and $\hLL(R)$ are isomorphic to each other for $q\not=\pm1$.

\begin{remark}\rm
Note, that for the $R$-matrix coming from the quantum group $U_q(sl(N))$ (the Drin\-feld-Jim\-bo $R$-matrix) the
defining relations (\ref{mREA}) turn into those (\ref{odin}) at the limit $q\rightarrow 0$ and the correponding modified
RE algebra $\hLL(R)$ tends to $U(gl(N))$.
\end{remark}

Also, we define the copies of the matrix $L$ ``shifted'' to higher positions\footnote{As we noticed above, the value of $p$ is not
specified to a concrete number.}:
$$
L_{\underline1}=L_{\overline 1}=L\otimes I^{\otimes (p-1)},\qquad L_{\overline k}=R_{k-1} L_{\overline{k-1}}R^{-1}_{k-1},\qquad
L_{\underline k}=R^{-1}_{k-1}L_{\underline{k-1}}R_{k-1},\qquad  k\ge 2.
$$
Moreover, we introduce a compact notation for products of the above coplies:
$$
L_{\,\overline{k\rightarrow r}} = L_{\,\overline k}L_{\,\overline{k+1}}\dots L_{\,\overline {r-1}}L_{\,\overline r}, \qquad
L_{\,\un{r\rightarrow k}} = L_{\,\un r}L_{\,\un{r-1}}\dots L_{\,\un {k+1}}L_{\,\un k},\quad k<r.
$$

Now, following \cite{IP} we define a map $ch:\,\H_n(q)\rightarrow \LL(R)$ $\forall n\ge 1$:
\be
\forall\, z\in \H_n(q):\quad  z\mapsto   ch(z)=\langle \rho_R(z)\, L_{\overline{1\rightarrow n}}\rangle_{1\dots n}.
\label{char-map}
\ee
As was proved in \cite{IP}, the family of all such images for all $n\geq 1$ is an associative subalgebra of the center of $\LL(R)$.
We we call it the {\em characteristic subalgebra}.

If in (\ref{char-map}) we take $z=e^\la_i$, the corresponding central element $s_\lambda(L) = ch(e^\lambda_i)$ does not
depend on the number $i$ and is called the Schur polynomial \cite{GPS2}. If we put $z=\tau_{n-1}...\tau_1$ or
$z= \tau_{n-1}^{-1}\dots \tau_1^{-1}$ we get respectively the power sums of first kind $p_n(L)$ or of the second kind $t_n(L)$.
In accordance with definition (\ref{char-map}) these central elements read:
\begin{eqnarray}
p_n(L)& =& \langle R_{n-1}\dots R_1 L_{\,\ov{1\to n}}\,\rangle_{1\dots n}=\langle R_{n-1}\dots R_1 L_{\,\un{n\to 1}}\,\rangle_{1\dots n},
\label{pk}\\
\rule{0pt}{5mm}
t_n(L) &=& \langle R^{-1}_{n-1}\dots R^{-1}_1 L_{\,\ov{1\to n}}\,\rangle_{1\dots n}=\langle R^{-1}_{n-1}\dots R^{-1}_1 L_{\,\un{n\to 1}} \,\rangle_{1\dots n}.\label{tk}
\end{eqnarray}
The last equalities  at this formulae are valid due to the following identity for generating matrices of RE algebra $\LL(R)$
(see \cite{IP}):
$$
L_{\,\overline{1\rightarrow k}} = L_{\,\underline{k\rightarrow 1}},\qquad \forall\,k\ge 2.
$$

Now, we assume all  symmetries $R$, we are dealing with, to be even of bi-rank $(m|0)$ (see \cite{GPS1}).  Consequently, the
generating matrix of the algebra $\LL(R)$  is subject to the following version of the  Cayley-Hamilton identity (see \cite{GPS3}):
\be
L^m+ \sum_{k=1}^m (-q )^ka_k(L)\,L^{m-k} = 0,
\label{q-CH}
\ee
where $a_k(L)$ is a $k$-th order elementary symmetric polynomial, i.e. the Schur polynomial, corresponding to the partition
$\la=(1,1,\dots,1)=(1^k)$.

Since the coefficients $a_k(L)$ of the polynomial (\ref{q-CH}) are central, it is natural,
similarly to the quantities $\hmu_i$ defined in introduction, to introduce new central elements $\{\mu_i\}_{1\le i\le m}$ as
solutions to the set of polynomial equations \cite{GPS5}:
\be
 qa_1(L) = \sum_{i=1}^m\, \mu_i, \qquad  q^2 a_2(L)  = \sum_{i<j}^m\, \mu_i\, \mu_j ,\qquad \dots \qquad  q^ma_m(L) =  \prod_i \mu_i.
\label{a-spec}
\ee
So, the elements $\mu_i$, called the {\it spectral values} (or eigenvalues) of the quantum matrix $L$, belong to an extended algebra
$\LL(R)[\mu_1,\dots \mu_m]$.

\begin{proposition}\label{p:2} \rm
The power sums of both kinds and elementary symmetric polynomials meet  the following identities for $\forall\, n\ge 1$:
\begin{eqnarray}
&&n_q a_n + \sum\limits_{k=1}^{n} (-1)^{k}q^{n-k} a_{n-k}\,p_k = 0\nonumber\\
&&n_q a_n + \sum\limits_{k=1}^{n} (-1)^{k} q^{k - n} a_{n-k}\, t_k = 0\nonumber\\
&&p_n = t_n + (q-q^{-1}) \sum\limits_{k = 1}^{n-1} p_{n-k} \,t_{k}.
\label{p-t}
\end{eqnarray}
In the above identities one should take into account that for Hecke symmetry with bi-rank $(m|0)$ the elementary
symmetric polynomials $a_n(L)\equiv 0$ if $n>m$.
\end{proposition}

\smallskip

\noindent
{\bf Proof.} The first and the second set of identities are called the quantum Newton relations, the proof of the first one can be found in \cite{GPS3}, the second set of identities can be proved in a similar way. We prove the third set of identities.

For this purpose we use the Hecke condition $R = R^{-1} + \omega I^{\otimes 2}$  (hereafter we use the notation $\omega = q-q^{-1}$)
and transform step by step the $R$ matrices in the expression for $p_k(L)$.
Thus, we get the following:
\begin{eqnarray*}
p_n&=&\langle L_{\overline{1 \rightarrow n}} R_{n-1} \dots R_1 \rangle_{1 \dots n} = \langle L_{\overline{1 \rightarrow n}}
R_{n-1}^{-1}R_{n-2} \dots R_1 \rangle_{1 \dots n} +
\omega \langle L_{\overline{1 \rightarrow n-1}} R_{n-2} \dots R_1  L_{\overline{n}}\rangle_{1 \dots n}\\
&=&\langle L_{\overline{1 \rightarrow n}} R_{n-1}^{-1}R_{n-2} \dots R_1 \rangle_{1 \dots n} + \omega \,p_{n-1} \,t_1 =
\langle L_{\overline{1 \rightarrow n}} R_{n-1}^{-1} R_{n-2}^{-1} R_{n-3} \dots R_1 \rangle_{1 \dots n}\\
&+&\omega \langle L_{\overline{1 \rightarrow n-2}} R_{n-3} \dots R_1 \rangle_{1 \dots n-2} \langle L_{\overline{n-1}}  L_{\overline{n}}
R_{n-1}^{-1} \rangle_{n-1\, n} + \omega \,p_{n-1} \,t_1\\
&=&\langle L_{\overline{1 \rightarrow n}} R_{n-1}^{-1} R_{n-2}^{-1} R_{n-3} \dots R_1 \rangle_{1 \dots n} + \omega ( \,p_{n-2}\, t_2
+ p_{n-1} \,t_1) =\dots\\
& =&  t_n + \omega \sum\limits_{k = 1}^{n-1} p_{n-k} t_{k}.
\end{eqnarray*}
Note that the system of relations (\ref{p-t}) is triangular (as well as the two other systems) and the elements $t_n$ can be expressed
as polynomials in $p_k$, $1\le k\le n$, by means of the Cramer's rule (and vise versa):
$$
t_n = (-1)^{n-1}
\det \left(\!\begin{array}{ccccc}
	p_1&1&0&...&0\\
	p_2&\omega p_1&1&...&0\\
	p_3&\omega p_2&\omega p_1&...&0\\
	...&...&...&...&1\\
	p_n&\omega e_{n-1}&\omega p_{n-2}&...&\omega p_1
\end{array}\!\right).
$$
This reflects the fact that finite sets $\{a_k(L)\}_{1\le k\le m}$, $\{p_k(L)\}_{1\le k\le m}$ and $\{t_k(L)\}_{1\le k\le m}$ are
different generating sets of the characteristic subalgebra. \hfill \rule{6.5pt}{6.5pt}

\smallskip

The spectral parameterization (\ref{a-spec}) results in the following parameterization of power sums (see \cite{GPS2}):
\be
p_n(L) =p_n(\mu) = \sum_{i=1}^m \mu_i^n\, d_i(\mu),\qquad d_i(\mu)=q^{-1}\prod_{j\not= i}\frac{\mu_i-q^{-2}\, \mu_j}{\mu_i-\mu_j}.
\label{pn-spec}
\ee

So, being expressed in terms of the quantum spectrum, the power sum $p_n(\mu)$ (\ref{pn-spec}) is proportional to the Hall-Littlewood
polynomial $q_n(\mu;t)$ \cite{M}, provided $t=q^{-2}$:
$$
q_n(\mu;t)_{\rule{0.4pt}{4mm}_{\,t=q^{-2}}}= \omega p_n(\mu).
$$

Introduce now the generating functions in formal variable $x$:
$$
A(x) =1+ \sum\limits_{k = 1}^{m} a_k x^k, \quad P(x) = \sum\limits_{k = 1}^{\infty} p_k x^k, \quad T(x) = \sum\limits_{k = 1}^{\infty} t_k x^k.
$$
Then  Proposition \ref{p:2} entails the following claim.
\begin{corollary}\rm The generating functions are connected with each other by the relations:
    \begin{equation}\label{AP}
    \omega A(qx) P(-x) = A\left(\frac{x}{q}\right) - A(qx),
    \end{equation}
    \begin{equation} \label{AT}
   \omega A\left(\frac{x}{q}\right) T(-x) = A\left(\frac{x}{q}\right) - A(qx),
    \end{equation}
\begin{equation}\label{PT}
	P(x) = T(x) + \omega P(x) T(x).
\end{equation}
\end{corollary}

Next, we use the spectral parameterization (\ref{a-spec}) in order to find the dependence of generating functions on the spectral
values $\mu_i$.

\begin{proposition} \rm The generating functions $P(x)$ and $T(x)$, being expressed via the eigenvalues $\mu_i$, read:
\be
1 + \om P(x) = \prod\limits_{i = 1}^{m} \frac{1 - q^{-2} \mu_i x}{1 - \mu_i x}, \qquad
1 - \om T(x) = \prod\limits_{i = 1}^{m} \frac{1 -  \mu_i x}{1 - q^{-2} \mu_i x}.
\label{firstsec}
\ee
\end{proposition}

\smallskip

\noindent
{\bf Proof.} With the use of parameterization (\ref{a-spec}) for the elementary symmetric functions $a_i(L)$ we obtain:
$$
A(qx) = 1+\sum\limits_{k = 1}^{m} x^k \Big(\sum\limits_{1 \leq i_1 <  \dots < i_k \leq m } \mu_{i_1} \dots \mu_{i_k}\Big)=
\prod\limits_{i=1}^{m} (1 + \mu_i x).
$$
Then the first formula from (\ref{firstsec}) is an immediate consequence of (\ref{AP}):
$$
1 + \om P(x) = \frac{A(-x /q)}{A(-qx)} = \prod\limits_{i = 1}^{m} \frac{1 - q^{-2} \mu_i x}{1 - \mu_i x}.
$$
The second formula in (\ref{firstsec}) follows from the relation $T(x) = P(x)(1 + \om P(x))^{-1}$ which is
just a reformulation of (\ref{PT}).\hfill\rule{6.5pt}{6.5pt}

\smallskip

It is not difficult to deduce from  formulas (\ref{firstsec}) and (\ref{pn-spec}) an expression for the power sums of the second
kind in terms of spectral values:
\be
t_n(L) = t_n(\mu) = \sum\limits_{i=1}^m (q^{-2} \mu_i)^n \widetilde{d}_i,\qquad  \widetilde{d}_i =q \prod\limits_{j \neq i}
\frac{\mu_i - q^2 \mu_j}{\mu_i - \mu_j}.
\label{tn-spec}
\ee

Now, we consider some elements of the representation theory of the algebras $\LL(R)$. To formulate necessary results
we introduce some convenient matrix notation.

First, the set of basis vectors of the space $V^{\otimes k}$, $k\ge 1$, will be denoted as
$$
\{x_i\otimes x_j\otimes \dots \otimes x_s\}\quad \Leftrightarrow \quad \x_1\x_2\dots \x_{k-1}\x_k.
$$
That is, the symbol $\x_i$ stands for an arbitrary basis vector of $V$  located at the $i$-th position in the tensor basis of $V^{\otimes k}$.

Moreover, we use the matrix summation rule: if an expression contains objects with the same numbers of matrix spaces, then
the summation over the corresponding indices is assumed. For example, the action of the operator
$R_{12} = R\otimes \mathrm{Id}_V$ on an arbitrary basis vector $x_i\otimes x_j\otimes x_k$ of the space $V^{\otimes 3}$
is described by the following matrix:
$$
\sum_{a,b = 1}^N R_{\,ij}^{\,\,ab}\,x_a\otimes x_b\otimes x_k\quad \Leftrightarrow\quad R_{12}\,\x_1\x_2\x_3.
$$

In addition to the space $V$, we consider its dual space $V^*$ and fix the right dual basis $\{y^i\}_{1\le i\le N}$ of $V^*$:
$$
\Omega(x_i, y^j) = \delta_i^j,
$$
where $\Omega: V\otimes V^*\rightarrow {\Bbb C}$ is a fixed non-degenerate bilinear form.

In the same way as above we introduce the symbols $\y^i$ for an arbitrary basis vectors at the $i$-th position of the tensor
product $V^{*\otimes k}$. Here we should emphasize that we accept reverse enumeration (from the right to the left)
of positions in the tensor product $V^{*\otimes k}$ of the dual spaces:
$$
\{y^i\otimes y^j\otimes \dots\otimes y^s\}\quad \Leftrightarrow \quad \y^k\y^{k-1}\dots \y^2\y^1.
$$
In accordance with this agreement, the action of the operator $R_{12}$ on an arbitrary basis vector $y^i\otimes y^j\otimes y^k$
of the space $V^{*\otimes 3}$ is described by the following matrix:
$$
\sum_{a,b = 1}^N y^i\otimes y^a\otimes y^bR_{\,ba}^{\,kj}\quad \Leftrightarrow\quad \y^3\y^2\y^1R_{12}.
$$

Now we can formulate the basic results on the representation of the RE algebra $\LL(R)$ in the spaces $V$ and $V^*$
(see \cite{GPS1,GPS4} for the proofs and detail).

\smallskip

\begin{proposition} \rm \phantom{a}

\noindent
{\bf 1. }The linear actions of the generators $l_i^j$ of the algebra $\LL(R)$ on the basis vectors of the spaces $V$ and $V^*$
in accordance with the formulas
\be
L_1R_1\triangleright \x_1= R^{-1}_1  \x_1\quad\Leftrightarrow\quad L_{\un 2} \triangleright \x_1 =J_2^{-1} \x_1,
\label{obt}
\ee
$$
L_2\triangleright \y^1= \y^1 J_2.
$$
define representations of $\LL(R)$. Here the symbol $\triangleright$ denotes the action of a linear operator,
$J_2 = R_{12}^2$.

\smallskip

\noindent
{\bf 2.} In the spaces $V^{\otimes k}$ and $V^{*\otimes k}$ the representations of $\LL(R)$ are defined by the following
linear actions of generators:
\be
L_{\un{k+1}}\triangleright \x_1 \x_2 \dots   \x_k=J_{k+1}^{-1}\,\x_1\x_2  \dots   \x_k.
\label{L-V-rep}
\ee
\be
L_{k+1} \triangleright \y^k\y^{k-1} \dots  \y^1 = \y^k\y^{k-1} \dots  \y^1 J_{k+1} .
\label{L-V-du}
\ee
\end{proposition}

\smallskip

To present the formula of action of arbitrary monomials in RE algebra generators on basis vectors of $V^{\otimes k}$
and $V^{*\otimes k}$ we introduce some more matrix objects.

Let us define a matrix $J_p^{\,\uparrow k}$ of the Jucys-Murphy element $j_p$ ``shifted'' by $k$ positions:
$$
J_p^{\,\uparrow k} = R_{k+p-1}\dots R_{k+2}R_{k+1}^2R_{k+2}\dots R_{k+p-1},\qquad \forall\,p\ge 2, \,\,\forall\,k\ge 0.
$$
One can easily show that matrices $J_p^{\,\uparrow k}$ with the same fixed shift $k$ and different $p$ commute with each other.

Now, using the action (\ref{L-V-rep}) of a single generator of $\LL(R)$ we can explicitly write down the action of any
$n$-th order monomial in generators:
\be
L_{\,\underline{k+n\rightarrow k+1}}\triangleright \x_1\x_2\dots \x_k =
\prod_{a=1}^nJ_{k+a}^{-1}\prod_{b=2}^nJ_b^{\,\uparrow k} \x_1\x_2\dots \x_k.
\label{mon-act}
\ee
The braid relation on $R$-matrices (\ref{braid-i}) allows one to prove the following commutativity:
$$
J_{k+a}^{-1}J_b^{\,\uparrow k} = J_b^{\,\uparrow k}J_{k+a}^{-1}, \quad \forall\, a>b.
$$
Using this commutativity we can rewrite the matrix products in the right hand side of (\ref{mon-act}):
\be
\prod_{a=1}^nJ_{k+a}^{-1}\prod_{b=2}^nJ_b^{\,\uparrow k}  = J_{k+1}^{-1}\left(J_{k+2}^{-1}J_2^{\,\uparrow k}\right)
\left(J_{k+3}^{-1}J_3^{\,\uparrow k}\right)\dots \left(J_{k+n}^{-1}J_n^{\,\uparrow k}\right) :=
J_{k+1}^{-1}\prod_{a=2}^{n}\left(J_{k+a}^{-1}J_a^{\,\uparrow k}\right).
\label{prod-re}
\ee
In the last formula we should always keep in mind the ordering in the product of noncommutative matrices:
the index $a$ increases from the left to the right.

To generalize (\ref{L-V-du}) to the action of higher order monomials  we define ``shifted'' copies
of the generating matrix:
\be
L_{\underline n}^{\,\uparrow k} = R_{k+n-1}^{-1}\dots R_{k+2}^{-1}R_{k+1}^{-1}L_{k+1}R_{k+1}R_{k+2}\dots
R_{k+n-1},\quad \forall\,n\ge 1, \,\,\forall\, k\ge 0.
\label{Lcop-sh}
\ee
In an similar way one can define the shifted copies $L_{\overline n}^{\,\uparrow k}$.

Now, by a straightforward calculation with multiple application of (\ref{L-V-du}) we find an analog of (\ref{mon-act})
but for basis vectors of $V^{*\otimes k}$
\be
L_{\,\underline{n\rightarrow 1}}^{\,\uparrow k}\triangleright \y^k\dots \y^1 = \y^k\dots \y^1 \prod_{b=2}^n(J_b^{-1})^{\,\uparrow k}
\prod_{a=1}^nJ_{k+a}  = \y^k\dots \y^1\prod_{a=2}^n\left((J_a^{-1})^{\,\uparrow k}J_{k+a}\right) J_{k+1}.
\label{mon-act*}
\ee
Note, that in the last expression the product of noncommutative matrices $(J_a^{-1})^{\,\uparrow k}J_{k+a}$ is
organized in reverse order comparing with (\ref{prod-re}): the index $a$ {\it decreases} form the left to the right.

The $\LL(R)$-modules $V^{\otimes k}$ and $V^{*\otimes k}$ for $k\ge 2$ are reducible and can be decomposed
into direct sums of (conjecturally\footnote{For the Drinfeld-Jimbo $R$-matrix defining the structure of the quantum group
$U_q(sl(N))$ the $\LL(R)$-modules $V^\lambda_i$ are irreducible. For an arbitrary skew-invertible Hecke symmetry $R$
the irreducibility of $V^\lambda_i$ is an open question.}) irreducible $\LL(R)$-modules.
These modules are images of orthogonal projectors obtained by $R$-matrix representation of the primitive idempotents
of the Hecke algebra ${\cal H}_k(q)$  (see \cite{GPS1} for detail):
\begin{equation}
V^{\otimes k} = \bigoplus_{\lambda\vdash k}\bigoplus_{i=1}^{d_\lambda}V^\lambda_i, \quad
V^\lambda_i = E^\lambda_i(R)\triangleright V^{\otimes k} = \mathrm{Im}(\rho_R(e^\lambda_{ii})),
\label{Vk-decom}
\end{equation}
where $e^\lambda_{ii}$ is a primitive idempotent of the Hecke algebra $\H_k(q)$, $E^\la_i(R)=\rho_R(e^\la_{ii})$ is a polynomial
in operators $R_i$, $1\le i\le k-1$. Note that the spaces $V^\lambda_i$ with the same $\lambda$ and different $i$ are equivalent
as $\LL(R)$-modules.

Analogous decomposition takes place for $V^{*\otimes k}$
\be
V^{*\otimes k} = \bigoplus_{\lambda\vdash k}\bigoplus_{i=1}^{d_\lambda}V^{*\lambda}_i, \quad
V^{*\lambda}_i  = \mathrm{Im}(\rho_R(e^\lambda_{ii}))
\label{Vdu-decom}
\ee

Though we are  not able to prove the irreducibility of $\LL(R)$-modules $V^\lambda_i$ and $V^{*\lambda}_i$, our computations
show that any element $a$ of the central characteristic subalgebra of $\LL(R)[\mu_1,\dots,\mu_m]$ becomes a
scalar operator being represented in modules $V^\la_i$ with the use of formulas (\ref{L-V-rep}):
$$
a\mapsto \chi_\la(a)\,\mathrm{Id}_{V^\lambda_i}.
$$
The factor $\chi_\la(a)$ will be called the character of the element $a$. The characters do not depend on the index $i$ at fixed $\lambda$
as they should do for the equivalent representations.

Upon representing the algebra $\LL(R)$ in the modules $V_\la^*$, we get characters $\chi_\la^*(a)$.

\begin{remark} \rm
For the Hecke symmetry $R$ coming from the quantum group $U_q(sl(N))$ it is possible to construct a homomorphism
of the algebra $\LL(R)$ into $U_q(sl(N))$ by using the method of \cite{RTF}. Thus, any $U_q(sl(N))$-module automatically
acquires the structure of an $\LL(R)$-module. Namely, this approach was used in \cite{JLM} in order to compute the characters
of the power sums.
\end{remark}

\section{Characters of spectral values}

In this section we present our main result: we explicitly write down the characters $\chi_\lambda(\mu_i)$ and
$\chi^*_\lambda(\mu_i)$ of central
spectral values of the quantum matrix in modules $V^\lambda_i$ and $V^{*\lambda}_i$. Since any element of the central
characteristic subalgebra can be written as a symmetric polynomial in spectral values $\mu_i$ we get a way to calculate
characters of any such an element in each module $V^\lambda_i$. In addition, our approach is universal in the sense that
it works for any Hecke symmetry $R$.

\smallskip

\begin{theorem}\label{th:7}\rm
The central spectral values $\{\mu_i\}_{1\le i\le m}$ defined in (\ref{a-spec}) are represented by scalar operators in $\LL(R)$-modules
$V^\lambda_a\subset V^{\otimes k}$ and $V^{*\lambda}_a\subset V^{*\otimes k}$, $\lambda\vdash k$:
\be
\mu_i\,\triangleright_{\,\rule{0.4pt}{4mm}{\,V^\lambda_a}} = \chi_\lambda(\mu_i)\,\mathrm{Id}_{V^\lambda_a}, \qquad
\mu_i\,\triangleright_{\,\rule{0.4pt}{4mm}{\,V^{*\lambda}_a}} = \chi^*_\lambda(\mu_i)\,\mathrm{Id}_{V^{*\lambda}_a}.
\label{mu-restr}
\ee
The characters do not depend on the index $a$ and are completely determined by the components of the partition
$\lambda= (\lambda_1\ge \lambda_2\ge\dots\ge \lambda_m\ge 0)$:
\begin{equation}
\chi_\lambda(\mu_i)= q^{-2(\lambda_i-i+m)}, \qquad \chi^*_\lambda(\mu_i) = q^{2({\lambda}_i - i + 1)}, \quad 1\le i\le m.
\label{chi-value}
\end{equation}
The form (\ref{mu-restr}) of restriction of the linear operators $\mu_i\triangleright$ is compatible with the actions
(\ref{L-V-rep}) and (\ref{L-V-du}) of the RE algebra on $V^{\otimes k}$ and $V^{*\otimes k}$.
\end{theorem}

\smallskip

To prove Theorem \ref{th:7} we need a number of technical results which we formulate and prove in the series of lemmas
below.

Introduce an axillary rational function in a formal variable $z$:
\be
S(z) =  \frac{z}{(1-z)^2}.
\label{S-ser}
\ee
Below we work with the generating functions $P(x)$ and $T(x)$ expressed in terms of the eigenvalues $\mu_i$ by means of 
formulas (\ref{firstsec}). For our purposes it is convenient to use a more detailed notation $P(x;\mu_1,\dots ,\mu_m)$ and
$T(x;\mu_1,\dots,\mu_m)$ for these generating functions.

\smallskip

\begin{lemma}\label{p:8} \rm
The following relations hold true:
 \be
 P(x; \mu_1, \dots,\mu_{k} q^{-2}, \dots, \mu_m) + \om S(q^{-2}x\mu_{k}) = P(x; \mu_1, \dots, \mu_m) \cdot
 (1 - \om^2 S(q^{-2}x\mu_{k}))
 \label{GenFuncP}
 \ee
 \be  T(x; \mu_1, \dots, \mu_{k} q^2, \dots, \mu_m) - \om S(x\mu_{k}) = T(x; \mu_1, \dots, \mu_m)  \cdot
 (1 - \om^2 S(x\mu_{k})).
 \label{GenFuncT}
 \ee
Here $\mu_{k}$ is an arbitrary fixed element from the set of eigenvalues $\{\mu_1,\dots.\mu_m\}$.
 \end{lemma}

\smallskip

\noindent
{\bf Proof.} The relations in the statement are proved in a similar way, so we check the second one.
Let us simplify the notations by setting
$$
T = T(x; \mu_1, \dots, \mu_m), \qquad T' = T(x; \mu_1, \dots, \mu_{k} q^2, \dots, \mu_m),  \qquad S = S(x\mu_{k})
=\frac{x\mu_k}{(1-x\mu_k)^2}.
$$
Then relation (\ref{GenFuncT}) is equivalent to the following
$$
1 - \om T'= (1 - \om T) (1 - \om^2 S).
$$
Next, it directly follows from (\ref{firstsec}) that:
$$
1 - \om T' = (1 - \om T) \frac{(1 - q^{-2} \mu_{k}x)(1 - q^2 \mu_{k}x)}{(1 -  \mu_{k}x)^2}.
$$
So, it suffices to check the relation
$$
1 - \om^2 S = \frac{(1 - q^{-2} \mu_{k}x)(1 - q^2 \mu_{k}x)}{(1 -  \mu_{k}x)^2},
$$
which is a straightforward consequence of the definition (\ref{S-ser}). \hfill\rule{6.5pt}{6.5pt}

\smallskip

Now we elaborate in some detail the RE algebra representation in spaces $V^{\otimes k}$ and $V^{*\otimes k}$ defined by
(\ref{L-V-rep}) and (\ref{L-V-du}) respectively. For this purpose we introduce the $N^k\times N^k$ matrices
$T_k^{\,n}$ and $P_k^{\,n}$:
\be
q^{2m (n-1)}T_k^{\,n} = \mathrm{Tr}_{R(k+1)}(J_{k+1}^{-n}) = \langle J_{k+1}^{-n}\rangle_{k+1},
\qquad P_k^{\,n} = \mathrm{Tr}_{R(k+1)} (J_{k+1}^n) = \langle J_{k+1}^n\rangle_{k+1}.
\label{matrix}
\ee
Recall that $R$-trace and the above bracket notation are defined in Section 2.

\smallskip

\begin{lemma}\label{lem:9}\rm
Let $n\ge 1$ be a fixed integer. Consider the first and the second kind power sums $p_n(L)$ and $t_n(L)$
defined in (\ref{pk}) and (\ref{tk}) respectively.

In the representation (\ref{L-V-rep}) the matrix
of the operator $t_n\triangleright\in\mathrm{End}(V^{\otimes k})$ in the tensor basis of $V^{\otimes k}$ is $T_k^{\,n}$:
$$
t_n\triangleright \x_1\x_2\dots \x_k = T_k^{\,n}\,\x_1\x_2\dots \x_k = q^{-2m(n-1)}\langle J_{k+1}^{-n}\rangle_{k+1}
\x_1\x_2\dots \x_k.
$$

In the representation (\ref{L-V-du}) the matrix of the operator $p_n\triangleright\in\mathrm{End}(V^{*\otimes k})$ in the tensor basis
of  $V^{*\otimes k}$ is $P_k^{\,n}$:
$$
p_n\triangleright \y^k\y^{k-1}\dots \y^1 = \y^k\y^{k-1}\dots \y^1P_k^{\,n} = \y^k\y^{k-1}\dots \y^1  \langle J_{k+1}^n\rangle_{k+1}.
$$
 \end{lemma}

\smallskip

\noindent
{\bf Proof.}  To calculate the matrix of the operator $t_n\triangleright$ in the tensor basis of the space $V^{\otimes k}$
we use the identity:
$$
t_n I_{1 \dots k} =  \langle L_{\underline{n \rightarrow 1}} R_{n-1}^{-1} \dots R_1^{-1} \rangle_{1 \dots n} I_{1 \dots k} =
\langle  L_{\underline{k+n \rightarrow k+1}} R_{k+n-1}^{-1} \dots R_{k+1}^{-1} \rangle_{k+1 \dots k+n}.
$$
Taking into account this identity, the action (\ref{mon-act}) and transformation (\ref{prod-re}) we find:
\begin{eqnarray*}
 t_n\!\!\!\!\!& \triangleright &\!\!\!\!\x_1 \dots \x_k = \langle L_{\underline{k+n \rightarrow k+1}}\triangleright \x_1 \dots \x_k
  R_{k+n-1}^{-1} \dots R_{k+1}^{-1} \rangle_{k+1 \dots k+n} \\
  &=&\!\!\!\!\langle J_{k+1}^{-1}\prod_{a=2}^{n}\left(J_{k+a}^{-1}J_a^{\,\uparrow k}\right)R_{k+n-1}^{-1}\dots R_{k+1}^{-1}
  \rangle_{k+n\dots k+1}\x_1\x_2\dots\x_k:=\Omega(R)\x_1\x_2\dots \x_k.
 \end{eqnarray*}
The matrix $\Omega(R)$ of the linear operator $t_n\triangleright$ is written in the last line as a multiple $R$-trace. We should prove
that $\Omega(R) = T_k^{\, n}$. To do so we explicitly calculate all $R$-traces except for that in the $(k+1)$-th position.

Note that in virtue of definitions of the Jucys-Murphy matrices we can simplify their products as follows:
\be
J_{k+a}^{-1}J_a^{\,\uparrow k} = R_{k+a-1}^{-1}\dots R_{k+1}^{-1}J_{k+1}^{-1}R_{k+1}\dots R_{k+a-1}.
\label{JJ}
\ee
First, we can calculate the $R$-trace in the $(k+n)$-th space. Appplying (\ref{JJ}) to the matrix $J_{k+n}^{-1}J_n^{\,\uparrow k}$,
we get:
\begin{eqnarray*}
\Omega(R) &:=&
\langle J_{k+1}^{-1}\prod_{a=2}^{n}\left(J_{k+a}^{-1}J_a^{\,\uparrow k}\right)R_{k+n-1}^{-1}\dots R_{k+1}^{-1}
\rangle_{k+n\dots k+1}\\
& =& \langle J_{k+1}^{-1}\prod_{a=2}^{n-1}\left(J_{k+a}^{-1}J_a^{\,\uparrow k}\right)R_{k+n-1}^{-1}\dots
R_{k+1}^{-1}J_{k+1}^{-1} \rangle_{k+n\dots k+1}\\
&=&q^{-2m}\langle J_{k+1}^{-1}\prod_{a=2}^{n-1}\left(J_{k+a}^{-1}J_a^{\,\uparrow k}\right)R_{k+n-2}^{-1}\dots
R_{k+1}^{-1}J_{k+1}^{-1} \rangle_{k+n-1\dots k+1}.
\end{eqnarray*}
In the third line of the above chain of transformations we calculate the $R$-trace in $(k+n)$-th space with the use of a general
formula:
\be
\langle X(R)_{1\dots p}R^{-1}_{p \,p+1}Y(R)_{1\dots p}\rangle_{p+1} = q^{-2m} X(R)_{1\dots p}Y(R)_{1\dots p}
\label{Rgen}
\ee
which is valid for any matrices $X(R)$ and $Y(R)$ provided they do not contain indices of $(p+1)$-th space.

Now, we apply (\ref{JJ}) to the next matrix $J_{k+n-1}^{-1}J_{n-1}^{\,\uparrow k}$, the chain of matrices
$R_{k+n-2}^{-1}\dots R_{k+1}^{-1}$ is reduced and we can calculate the $(k+n-1)$-th $R$-trace with the general formula
(\ref{Rgen}) thus coming to the intermediate result:
$$
\Omega(R) =
q^{-4m}\langle J_{k+1}^{-1}\prod_{a=2}^{n-2}\left(J_{k+a}^{-1}J_a^{\,\uparrow k}\right)R_{k+n-3}^{-1}\dots
R_{k+1}^{-1}J_{k+1}^{-2} \rangle_{k+n-2\dots k+1}.
$$
Repeating this procedure step by step, we calculate all $R$-traces except for $(k+1)$-th ending with the desired result:
$$
\Omega(R) = q^{-2m(n-1)}\langle J_{k+1}^{-n}\rangle_{k+1} = T_k^{\,n}.
$$

To calculate the matrix of the linear operators $p_n\triangleright$ we use the identity:
$$
p_n(L)I_{1\dots k} = \langle R_{k+1}R_{k+2}\dots R_{k+n-1}L_{\,\underline{n\rightarrow 1}}^{\,\uparrow k}\rangle_{k+n\dots k+1}.
$$

Taking into account this identity and the action (\ref{mon-act*}) we find:
\begin{eqnarray*}
	p_n\!\!\!\!\!& \triangleright &\!\!\!\!\y^k\dots \y^1 = \langle R_{k+1} \dots R_{k+n-1}
	L_{\underline{n \rightarrow 1}}^{\,\uparrow k}\triangleright \y^k\dots \y^1 \rangle_{k+1 \dots k+n} \\
	&=&\!\!\!\! \y^k\dots \y^1\langle R_{k+1}\dots R_{k+n-1} \prod_{a=2}^n\left((J_a^{-1})^{\,\uparrow k}J_{k+a}\right) J_{k+1}
	\rangle_{k+n\dots k+1} := \y^k\dots \y^1\tilde \Omega(R).
\end{eqnarray*}
Now we calculate the $R$-traces in the matrix $\tilde\Omega(R)$ in the same way as we did for the matrix $\Omega(R)$ above.
The only difference is that instead of the general formula (\ref{Rgen}) we use its analog
$$
\langle X(R)_{1\dots p}R_{p \,p+1}Y(R)_{1\dots p}\rangle_{p+1} = X(R)_{1\dots p}Y(R)_{1\dots p}
\label{genR}
$$
which does not contain the multiplier $q^{-2m}$.

As a result we get
$$
\tilde\Omega(R) = \langle J_{k+1}^n\rangle_{k+1} = P^{\,n}_k
$$
thus completing the proof of the lemma.   \hf

\smallskip
The last auxiliary technical result concerns the relations among matrices $P_k^{\,n}$ and $T_k^{\,n}$ with different values of k.

\begin{lemma}\label{l:10} \rm The following identities hold true:
\be
P_k^{\,n} = P_{k-1}^{\,n} + \om n J_k^{n} + \om^2\sum\limits_{s=1}^{n-1} s P^{\,n-s}_k J_k^{s},
\label{Pn-id}
\ee
\be
T_k^{\,n} = T_{k-1}^{\,n} - \om n (q^{2m} J_k )^{-n} + \om^2 \sum\limits_{s=1}^{n-1} s T_k^{\,n-s} (q^{2m} J_k )^{-s}.
\label{Tn-id}
\ee
\end{lemma}

\smallskip

\noindent
{\bf Proof.} Take into account the matrix equality:
$$
R_kJ_{k+1}= R_k^2J_kR_k = J_kR_k + \omega J_{k+1},
$$
where the Hecke condition $R_k^2 = I^{\otimes 2} + \om\, R_k$ was applied.

Using this equality as the base of induction in $n$, it is not difficult to prove the matrix identity:
\be
R_k J_{k+1}^n = J_k^{n} R_k + \om \sum\limits_{s=1}^{n} J_{k+1}^s J_{k}^{n-s}.
\label{RJn}
\ee

We apply the identity (\ref{RJn}) in order to transform the matrix power $J_{k+1}^n$:
\begin{eqnarray*}
J_{k+1}^n &=&  R_k J_k (R_k J_{k+1}^{n-1}) =  R_k J_k^{n} R_k + \omega \sum\limits_{s=1}^{n-1}R_k J_{k+1}^s J_k^{n-s}\\
&=& R_kJ_k^nR_k  + \omega \sum\limits_{s=1}^{n-1} J_{k}^{s} R_k J_{k}^{n-s} + \omega^2 \sum_{s =1}^{n-1}\sum_{r=1}^s
J_{k}^{n-r} J_{k+1}^{r}\\
&=& R_kJ_k^nR_k  + \omega \sum\limits_{s=1}^{n-1} J_{k}^{s} R_k J_{k}^{n-s}+\omega^2\sum_{s=1}^{n-1}s J_k^sJ_{k+1}^{n-s}.
\end{eqnarray*}
In the last line of the above transformations we use the Hecke condition in the form $R_k = R_k^{-1}+\omega I^{\otimes 2}$
in order to rewrite the frist term:
$$
R_k J_k^{n} R_k =  R_k J_k^n R_k^{-1}+\omega R_k J_k^n
$$
Then, applying the $R$-trace in the $(k+1)$-th space we get (\ref{Pn-id}). In so doing, for the first term we use an important
property of the $R$-trace:
$$
\langle R_kX(R)_{1\dots k}R_k^{-1})\rangle_{k+1} = \langle X(R)_{1\dots k}\rangle_k
$$
valid for any matrix $X(R)$ which does not contain indices of the $(k+1)$-th space.

The identity (\ref{Tn-id}) can be proved in the same way by transforming the matrix power $J_{k+1}^{-n}$ and applying
the $R$-trace in the $(k+1)$-th space at the end of calculations. \hfill\rule{6.5pt}{6.5pt}

\smallskip

\noindent
{\bf Proof of the Theorem \ref{th:7}.}

\smallskip

We are going to prove that the values of characters $\chi_\lambda(t_n)$ and $\chi^*_\lambda(p_n)$
corresponding to the restrictions of the linear operators $t_n\triangleright$ and $p_n\triangleright$ on $\LL(R)$-modules
$V^\lambda_a$ and $V^{*\lambda}_a$ respectively read as follows:
\be
t_n\,\triangleright_{\,\rule{0.4pt}{4mm}{\,V^\lambda_a}} = \chi_\lambda(t_n)\mathrm{Id}_{\,V^\lambda_a}, \qquad
\chi_\lambda(t_n) = t_n(\mu_1, \dots, \mu_m)_{\,\rule{0.4pt}{4mm}{\,\,\mu_i = q^{-2 ( {\lambda}_i - i + m)} }},
\label{TnAct}
\ee
\be
p_n\,\triangleright_{\,\rule{0.4pt}{4mm}{\,V^{*\lambda}_a}} = \chi_\lambda^*(p_n)\mathrm{Id}_{\,V^\lambda_a}, \qquad
\chi_\lambda^*(p_n) = p_n(\mu_1, \dots, \mu_m)_{\,\rule{0.4pt}{4mm}{\,\,\mu_i = q^{2 ( {\lambda}_i - i + 1)} }}.
\label{PnAct}
\ee
Here the polynomials $p_n(\mu)$ and $t_n(\mu)$ are defined in (\ref{pn-spec}) and (\ref{tn-spec}).

Note that the spectral parameterizations (\ref{pn-spec}) and (\ref{tn-spec}) of the power sums are in one-to-one
correspondence with the parametrerization (\ref{a-spec}). This is a consequence of the fact that the sets
$\{a_k(L)\}_{1\le k\le m}$ and $\{p_k(L)\}_{1\le k\le m}$ are the generating sets of the characteristic subalgebra
and, therefore, the values of power sums defines the values of $\mu_i$ up to an ordering. So,
if the relations (\ref{TnAct}) and (\ref{PnAct}) are indeed fulfilled, then it leads to the conclusion that above characters
of spectral values solve the defining system (\ref{a-spec}).

We will prove (\ref{TnAct}) and (\ref{PnAct}) by double induction in the parameter $k=\sum_i \la_i$ parameterizing
the tensor power of the representation space, and for each  $k$ we realize the internal induction in $n$ parameterizing the
order of the power sums $p_n(L)$ and $t_n(L)$.

According to Lemma \ref{lem:9} and definitions $(\ref{Vk-decom})$, $(\ref{Vdu-decom})$ the matrices of operators
$t_n\triangleright$ and $p_n\triangleright$ restricted to the $\LL(R)$-invariant subspaces $V^{\lambda}_a$ and
$V^{*\lambda}_a$ are respectively $T_k^{\,n}E_a^{\lambda}(R)$ and $P_k^{\,n}E_a^{\lambda}(R)$. Here
$E_{a}^{{\lambda}}(R) = \rho_R(e_{aa}^\lambda)$ is an idempotent, corresponding to the diagram $\la \vdash k$ (we
use the same term for the idempotent of the Hecke algebra and for its matrix image in the $R$-matrix representation).

Thus, as the step of induction, we assume that the matrix product of $P^n_{k-1}$ or $T^n_{k-1}$ and the idempotent
corresponding to the Young table with omitted box with number $k$ as well as the matrix products
$P^l_{k}E_{a}^{{\lambda}}(R)$ or $T^l_{k}E_{a}^{{\lambda}}(R)$ with $l<n$ have the necessary form (\ref{TnAct}) or
$(\ref{PnAct})$.

In the standard Young table corresponding to the idempotent $E_a^\lambda$ the largest number $k$ is placed at the
ending box of some row. Let $i_0$ be the number of this row. Then the content of the corresponding box is
$c(k) = {\lambda}_{i_0} - i_0$,
and consequently $J_k E_{a}^{{\lambda}} = q^{2({\lambda}_{i_0} - i_0)}E_{a}^{{\lambda}}$ (see \cite{OP} for detail).

Using this consideration and Lemma \ref{l:10}, we get
\be
T_k^n E_{a}^{{\lambda}} = T_{k-1}^{n} E_{a}^{ {\lambda}} - \om n  (q^{-2( {\lambda}_{i_0} - i_0 + m)})^{n} E_{a}^{ {\lambda}} +
\om^2 \sum\limits_{s=1}^{n-1} s (q^{-2( {\lambda}_{i_0} - i_0 + m)})^{s} T_k^{n-s} E_{a}^{ {\lambda}}
\label{TknE}
\ee
and
\be
P_k^n E_{a}^{{\lambda}} = P_{k-1}^{n} E_{a}^{{\lambda}} + \om n  (q^{2({\lambda}_{i_0} - i_0)})^{n} E_{a}^{{\lambda}} + \om^2 \sum\limits_{s=1}^{n-1} s (q^{2({\lambda}_{i_0} - i_0)})^{i} P_k^{n-s} E_{a}^{{\lambda}}.
\label{PknE}
\ee

Next, the relations (\ref{GenFuncP}), (\ref{GenFuncT}) on the generating functions lead to the following
identities for the power sums:
\be
t_n(\mu_1, \dots, \mu_m) = t_n(\mu_1, \dots, q^2 \mu_{i_0}, \dots, \mu_m) - \om n \mu_{i_0}^{n} + \om^2 \sum_{s=1}^{n-1}
s \mu_{i_0}^s t_{n-s}(\mu_1, \dots, \mu_m),
\label{t-mu-id}
\ee
and
\begin{eqnarray}
p_n(\mu_1, \dots, \mu_m) &=& p_n(\mu_1, \dots, q^{-2} \mu_{i_0}, \dots, \mu_m)\nonumber\\
& + &\om n (q^{-2}\mu_{i_0})^{n} + \om^2 \sum\limits_{s=1}^{n-1} s (q^{-2} \mu_{i_0})^s p_{n-s}(\mu_1, \dots, \mu_m).
\label{p-mu-id}
\end{eqnarray}

Compare now the formulas (\ref{TknE}) and (\ref{t-mu-id}). First, we point out that
the right hand side of (\ref{TknE}) equals to $E_a^\lambda(R)$ multiplied by a scalar factor. By definition this factor
is a character $\chi_\lambda(t_n)$ since at the left hand side of (\ref{TknE}) we have the matrix of $t_n\triangleright$
restricted on $V^\lambda_a$. Indeed, the matrices $T^{\,n-s}_kE_a^\lambda$ are of the form $\chi_\lambda(t_{n-s})$
by assumption of induction, the same is true for the first term $T_{k-1}^nE_a^\lambda(R)$. The only subtle point is that
the product of $T_{k-1}^{n}E_{a}^{ {\lambda}}(R)$ contains the character corresponding to the table with eliminated
box with number  $k$, which, in turn, equivalent to the change of $\mu_{i_0}$ to $q^2\mu_{i_0}$. But from the other hand,
having represented the central elements of algebraic identity (\ref{t-mu-id}) in the module $V^\lambda_a$, we get the same
connection among their characters. So, we conclude that the formula (\ref{TnAct}) indeed takes place.

For polynomials $p_n$ the consideration is the same, except for the change of $\mu$ corresponding to eliminated box with
the number $k$: one should substitute $q^{-2}\mu_{i_0}$ instead of $\mu_{i_0}$.

So, the induction step is proved and it remains to chek a few particular cases forming the base of induction.

If $n = 1$ we get the relation (see \cite{GPS1})
$$
T_k^1 E^{ {\lambda}}_{a} = \langle J_{k+1}^{-1} \rangle_{k+1} E^{ {\lambda}}_{a} = \frac{1}{q} \sum\limits_{s=1}^{m}
 q^{-2( {\lambda}_s - s + m)}E^{ {\lambda}}_{a}.
$$
Also we use the formula $P_k^1 = \langle J_{k+1} \rangle_{k+1} = \frac{m_q}{q^m} + \om \sum\limits_{i=1}^{k} J_i$.
Multipluing by the idempotent $E_a^\lambda(R)$ we find:
\begin{eqnarray*}
P_k^1E^{ {\lambda}}_{a}(R) & =& \Big(\,\frac{m_q}{q^m} + \om\sum\limits_{i=1}^{k} q^{2 c(i)}\Big)E^{ {\lambda}}_{a}(R)
 = \Big(\,\frac{m_q}{q^m} +  \om\sum\limits_{i=1}^{m} q^{2(1-i)} \sum\limits_{j=0}^{ {\lambda}_i -1} q^{2j}\Big)E^{ {\lambda}}_{a}(R) \\
 &=&\Big(\,\frac{m_q}{q^m} - \frac{1}{q} \sum\limits_{i=1}^{m} q^{2(1-i)} + \frac{1}{q} \sum\limits_{i=1}^{m} q^{2(\lambda_i-i+1)}\Big)
 E^\lambda_{a}(R) = \frac{1}{q} \sum\limits_{i=1}^{m} q^{2( {\lambda}_i - i + 1)}E^{ {\lambda}}_{a}(R).
 \end{eqnarray*}
Thus, we get a conclusion for $t_1$ and $p_1$, that is, the base of the internal induction in the parameter $n$ (the order of
the power sums) is established.

Now let us check the base of external induction, that is $k=0$ while $n$ is arbitrary. In this case $T_0^{\,n}$ is calculated easily:
$$
T_0^{\,n} = q^{-2m(n-1)}\langle J_1^n\rangle_1 = q^{m(1-2n)}m_q.
$$
Since $k=0$ then $\lambda = 0$ and the characters of the eigenvalues read $\chi_0(\mu_i) = q^{-2(m-i)}$.  Omitting
for simplicity the notation $\chi_0$ we see that $\mu_{i+1} = q^2\mu_i$, $1\le i\le m$. As a consequence, the multiplicities
$\tilde d_i$ in (\ref{tn-spec}) are all equal to zero except for $\tilde d_1$ and we get for $T_0^{\,n}$ the same result as above:
$$
\tilde d_1 = q\prod_{i=2}^m\frac{q^{2(1-m)} - q^{2(i+1 - m)}}{q^{2(1-m)} - q^{2(i - m)}} = q^mm_q\quad \Rightarrow\quad
T_0^{\,n} = (q^{-2}\mu_1)^n\tilde d_1 = q^{m(1-2n)}m_q.
$$

As for $P_0^{\,n}$ the verification of the base of induction is also straightforward. For $k = 0$ and an arbitrary $n$ the
characterss of spectral values $\chi_0^*(\mu_i)=q^{2(1-i)}$, so in virtue of analogous property $\mu_{i+1} = q^{-2}\mu_i$
all the multiplicities $d_i$ in (\ref{pn-spec}) are equal to zero except for $d_1$ and we get for $P_0^{\,n}$:
$$
P_0^{\,n} = (\mu_1)^nd_1 = q^{-1}\prod_{i=2}^m\frac{1 - q^{-2i}}{1 - q^{-2(i - 1)}} =
q^{-m}m_q = \langle J_1^n\rangle_1.
$$
This completes the verification of the base of induction and the proof of Theorem \ref{th:7}.
\hfill\rule{6.5pt}{6.5pt}

\smallskip

As an example of application of Theorem \ref{th:7} we present the characters $\chi_\lambda(p_n)$ and $\chi_\lambda^*(p_n)$.

\smallskip

\begin{corollary}\rm
On applying to the spectral parameterization (\ref{pn-spec}) the results of Theorem \ref{th:7} one gets the characters
of the power sums represented in $V^\lambda_a$ and $V^{*\lambda}_a$ with the actions of generators (\ref{L-V-rep})
and (\ref{L-V-du}):
\begin{equation}
\chi_\lambda(p_n) = q^{-m(2n+1)}\sum_{i=1}^mq^{-2n\ell_i}\prod_{j\not=i}\frac{(\ell_i-\ell_j-1)_q}{(\ell_i-\ell_j)_q},
\label{pn-chi}
\end{equation}
\begin{equation}
\chi^*_\lambda(p_n) = q^{2n-m}\sum_{i=1}^mq^{2n\ell_i}\prod_{j\not=i}\frac{(\ell_i-\ell_j+1)_q}{(\ell_i-\ell_j)_q},
\end{equation}
where $\ell_i = \lambda_i-i$ and a $q$-number is defined as $n_q = (q^n-q^{-n})/(q-q^{-1})$ $\forall\, n\in{\Bbb Z}$.
Note, that in \cite{JLM} the relation (\ref{pn-chi}) was obtained for the $U_q(gl(N))$ case (in this case $m=N$) by different methods.
\end{corollary}

In conclusion, we extract from our formulas the classical results obtained in \cite{PP}. At first, we have to pass to the generating
matrix $\hL$ of the corresponding  modified RE algebras $\hLL(R)$ by the shift of generators (\ref{isom}).
The matrix $\hL$ is also subject to a version of the CH identity. It can be readily obtained from (\ref{CH}) and (\ref{isom}).
Let $\hmu_i$ be the eigenvalues of the characteristic polynomial for the matrix $\hL$. They are related to the eigenvalues of the
matrix $L$ by the formula
$$
\mu_i=1-\omega \hmu_i.
$$

This relation enables us to compute the quantities $\chi_\la({\hat{\mu}_i})$. Namely, we have
$$
\chi_\la({\hat{\mu}_i})=\frac{1- q^{-2\, (\la_i+m-i)}}{q-\qq} = q^{-(\lambda_i+m-i)}(\lambda_i+m-i)_q,\quad 1\leq i \leq m.
$$
The multiplicities $d_i(\hat\mu)$ read as follows:
$$
d_i(\hat\mu)= q^{-1}\prod_{j\not=1}\frac{\hat \mu_i -\hat \mu_j - q^{-1}}{\hat \mu_i - \hat\mu_j}.
$$

If  $R\to P$ as $q\to 1$, then by passing to the limit $q\rightarrow 1$ we obtain  the quantities
$\chi_\la(\hat{\mu}_k)$ and $d_i(\chi_\lambda(\hat\mu))$for the matrix $\hL$ generating the algebra $U(gl(N))$ (in this case $m=N$):
$$
\chi_\la({\hat{\mu}_i})=\la_i-i+N,\quad d_i(\hat\mu) = \prod_{j\not=i}\frac{(\lambda_i -\lambda_j+j-i-1)}{(\lambda_i -\lambda_j+j-i)},
\quad 1\leq i\leq N.
$$
Thus, we recover the result of \cite{PP} for $\mathrm{Tr}\hL^n=\sum_i\hat\mu_i^nd_i(\hat\mu)$ (though, in the cited paper
the notion of the eigenvalues of the matrix $\hL$ was not used).

\end{document}